# APPROXIMATION AND LEARNING BY GREEDY ALGORITHMS[1]

By Andrew R. Barron, Albert Cohen,
Wolfgang Dahmen and Ronald A. DeVore

*Yale University, Université Pierre et Marie Curie 175,
RWTH Aachen and University of South Carolina*

We consider the problem of approximating a given element $f$ from a Hilbert space $\mathcal{H}$ by means of greedy algorithms and the application of such procedures to the regression problem in statistical learning theory. We improve on the existing theory of convergence rates for both the orthogonal greedy algorithm and the relaxed greedy algorithm, as well as for the forward stepwise projection algorithm. For all these algorithms, we prove convergence results for a variety of function classes and not simply those that are related to the convex hull of the dictionary. We then show how these bounds for convergence rates lead to a new theory for the performance of greedy algorithms in learning. In particular, we build upon the results in [*IEEE Trans. Inform. Theory* **42** (1996) 2118–2132] to construct learning algorithms based on greedy approximations which are universally consistent and provide provable convergence rates for large classes of functions. The use of greedy algorithms in the context of learning is very appealing since it greatly reduces the computational burden when compared with standard model selection using general dictionaries.

**1. Introduction.** We consider the problem of approximating a function $f$ from a Hilbert space $\mathcal{H}$ by a finite linear combination $\hat{f}$ of elements of a given dictionary $\mathcal{D} = (g)_{g \in \mathcal{D}}$. Here, by *dictionary*, we mean any family of functions from $\mathcal{H}$. In this paper, this problem is addressed in two different contexts:

Received February 2006; revised March 2007.
[1]Supported by the Office of Naval Research Contracts ONR-N00014-03-1-0051, ONR/DEPSCoR N00014-03-1-0675 and ONR/DEPSCoR N00014-00-1-0470; the Army Research Office Contract DAAD 19-02-1-0028; the AFOSR Contract UF/USAF F49620-03-1-0381; the French–German PROCOPE contract 11418YB; and NSF Grants DMS-02-21642 and DMS-02-00187.

(i) Deterministic approximation: $f$ is a known function in a Hilbert space $\mathcal{H}$. The approximation error is naturally measured by $\|f - \hat{f}\|$, where $\|\cdot\|$ is the corresponding norm generated by the inner product $\langle \cdot, \cdot \rangle$ on $\mathcal{H}$, that is, $\|g\|^2 := \|g\|_{\mathcal{H}}^2 := \langle g, g \rangle$.

(ii) Statistical learning: $f = f_\rho$, where $f_\rho(x) = E(y|x)$ is the regression function of an unknown distribution $\rho$ on $X \times Y$, with $x \in X$ and $y \in Y$ respectively representing the feature and output variables, from which we observe independent realizations $(z_i) = (x_i, y_i)$ for $i = 1, \ldots, n$. The approximation error is now measured in the Hilbertian norm $\|u\|^2 := E(|u(x)|^2)$.

In either of these situations, we may introduce the set $\Sigma_N$ of all possible linear combinations of elements of $\mathcal{D}$ with at most $N$ terms and define the best $N$-term approximation error $\sigma_N(f)$ as the infimum of $\|f - \hat{f}\|$ over all $\hat{f}$ of this type,

$$(1.1) \qquad \sigma_N(f) = \inf_{\#(\Lambda) \leq N} \inf_{(c_g)} \left\| f - \sum_{g \in \Lambda} c_g g \right\|.$$

In the case where $\mathcal{D}$ is an orthonormal basis, the minimum is attained by

$$(1.2) \qquad \hat{f} = \sum_{g \in \Lambda_N(f)} c_g g,$$

where $\Lambda_N(f)$ corresponds to the coordinates $c_g := \langle f, g \rangle$ which are the $N$-largest in absolute value. The approximation properties of this process are well understood; see, for example, [8]. In particular, one can easily check that the convergence rate $\|f - \hat{f}\|_{\mathcal{H}} \lesssim N^{-s}$ is equivalent to the property that the sequence $(c_g)_{g \in \mathcal{D}}$ belongs to the weak space $w\ell_p$ with $1/p = 1/2 + s$. We recall here that $w\ell_p$ is the space of sequences $(c_g)_{g \in \mathcal{D}}$ such that the quasi-norm $\|(c_g)\|_{w\ell_p}$, defined by

$$(1.3) \qquad \|(c_g)\|_{w\ell_p}^p := \sup_{\eta > 0} \eta^p \#(\{g; |c_g| \geq \eta\}),$$

is finite. Note that $\|(c_g)\|_{w\ell_p}^p \leq \sum_{g \in \mathcal{D}} |c_g|^p := \|(c_g)\|_{\ell_p}^p$, so that $\ell_p \subset w\ell_p$ (see, e.g., the survey [8] or standard books on functional analysis for a more expanded discussion of $\ell_p$ and weak $\ell_p$ spaces). Here, and later in this paper, we use the notation $A \lesssim B$ to mean that $A \leq CB$ for some absolute constant $C$ that does not depend on the parameters which define $A$ and $B$.

One of the motivations for utilizing general dictionaries rather than orthonormal systems is that in many applications such as signal processing and statistical estimation, it is not clear which orthonormal system, if any, is best for representing or approximating $f$. Thus, dictionaries which are unions of several bases or collections of general waveforms are preferred. Some well-known examples are the use of Gabor systems, curvelets and wavepackets in signal processing and neural networks in learning theory. Moreover, in statistical learning, orthonormality is meaningful in the norm $\|u\|^2 := E(|u(x)|^2)$, which is usually out of reach since it depends on the unknown distribution of the random variable $x$. Therefore, non-orthonormal systems cannot be avoided in this context.

When working with dictionaries $\mathcal{D}$ which are not orthonormal bases, the realization of a best $N$-term approximation is usually out of reach from a computational point of view since it would require minimizing $\|f - \hat{f}\|$ over all $\hat{f}$ in an infinite or huge number of $N$-dimensional subspaces. *Greedy algorithms* or matching pursuits aim to build "suboptimal yet good" $N$-term approximations through a greedy selection of elements $g_k$, $k = 1, 2, \ldots$, within the dictionary $\mathcal{D}$ and to do so with a more manageable number of computations.

1.1. *Greedy algorithms.* Greedy algorithms have been introduced in the context of statistical estimation. They have also been considered for applications in signal processing [1]. Their approximation properties have been explored in [4, 9, 14, 15, 18, 19] for general bounded dictionaries along with various applications. A recent survey of the approximation properties of such algorithms is given in [21].

There exist several versions of these algorithms. The four most commonly used are the *pure greedy*, the *orthogonal greedy*, the *relaxed greedy* and the *stepwise projection* algorithms, which we respectively denote by the acronyms PGA, OGA, RGA and SPA. We describe these algorithms in the deterministic setting. We shall assume here and later that the elements of the dictionary are normalized according to $\|g\| = 1$ for all $g \in \mathcal{D}$ unless it is explicitly stated otherwise.

All four of these algorithms begin by setting $f_0 := 0$. We then recursively define the approximant $f_k$ based on $f_{k-1}$ and its residual $r_{k-1} := f - f_{k-1}$.

In the PGA and the OGA, we select a member of the dictionary as

$$(1.4) \qquad g_k := \underset{g \in \mathcal{D}}{\operatorname{Arg\,max}} |\langle r_{k-1}, g \rangle|.$$

The new approximation is then defined as

$$(1.5) \qquad f_k := f_{k-1} + \langle r_{k-1}, g_k \rangle g_k$$

in the PGA and as

$$(1.6) \qquad f_k = P_k f$$

in the OGA, where $P_k$ is the orthogonal projection onto $V_k := \operatorname{Span}\{g_1, \ldots, g_k\}$. It should be noted that when $\mathcal{D}$ is an orthonormal basis, both algorithms coincide with the computation of the best $k$-term approximation.

In the RGA, the new approximation is defined as

$$f_k = \alpha_k f_{k-1} + \beta_k g_k, \tag{1.7}$$

where $(\alpha_k, \beta_k)$ are real numbers and $g_k$ is a member of the dictionary. There exist many possibilities for the choice of $(\alpha_k, \beta_k, g_k)$, the most greedy being to select them according to

$$(\alpha_k, \beta_k, g_k) := \operatorname*{Arg\,min}_{(\alpha,\beta,g)\in\mathbb{R}^2\times\mathcal{D}} \|f - \alpha f_{k-1} - \beta g\|. \tag{1.8}$$

Other choices specify one or several of these parameters, for example, by taking $g_k$ as in (1.4) or by setting in advance the values of $\alpha_k$ and $\beta_k$; see, for example, [14] and [4]. Note that the RGA coincides with the PGA when the parameter $\alpha_k$ is set to 1.

In the SPA, the approximation $f_k$ is defined by (1.6), as in the OGA, but the choice of $g_k$ is made so as to minimize over all $g \in \mathcal{D}$ the error between $f$ and its orthogonal projection onto $\operatorname{Span}\{g_1, \ldots, g_{k-1}, g\}$.

Note that from a computational point of view, the OGA and SPA are more expensive to implement since at each step, they require the evaluation of the orthogonal projection $P_k f$ (and, in the case of SPA, a renormalization). Such projection updates are computed preferably using Gram–Schmidt orthogonalization (e.g., via the QR algorithm) or by solving the normal equations

$$G_k a_k = b_k, \tag{1.9}$$

where $G_k := (\langle g_i, g_j \rangle)_{i,j=1,\ldots,k}$ is the Gramian matrix, $b_k := (\langle f, g_i \rangle)_{i=1,\ldots,k}$ and $a_k := (\alpha_j)_{j=1,\ldots,k}$ is the vector such that $f_k = \sum_{j=1}^k \alpha_j g_j$.

In order to describe the known results concerning the approximation properties of these algorithms, we introduce the class $\mathcal{L}_1 := \mathcal{L}_1(\mathcal{D})$ consisting of those functions $f$ which admit an expansion $f = \sum_{g \in \mathcal{D}} c_g g$, where the coefficient sequence $(c_g)$ is absolutely summable. We define the norm

$$\|f\|_{\mathcal{L}_1} := \inf\left\{\sum_{g \in \mathcal{D}} |c_g| : f = \sum_{g \in \mathcal{D}} c_g g\right\} \tag{1.10}$$

for this space. This norm may be thought of as an $\ell_1$ norm on the coefficients in the representation of the function $f$ by elements of the dictionary; we emphasize that it is not to be confused with the $L_1$ norm of $f$. An alternate, and closely related, way of defining the $\mathcal{L}_1$ norm is by the infimum of numbers $V$ for which $f/V$ is in the closure of the convex hull of $\mathcal{D} \cup (-\mathcal{D})$. This is known as the "variation" of $f$ with respect to $\mathcal{D}$ and was used in [16, 17], building on the earlier terminology in [3].

In the case where $\mathcal{D}$ is an orthonormal basis, we find that if $f \in \mathcal{L}_1$, then

$$\sigma_N(f) = \left( \sum_{g \notin \Lambda_n(f)} |c_g|^2 \right)^{1/2}$$

$$(1.11) \qquad \leq \left( \|f\|_{\mathcal{L}_1} \min_{g \in \Lambda_n(f)} |c_g| \right)^{1/2}$$

$$\leq \|f\|_{\mathcal{L}_1} N^{-1/2}.$$

For the PGA, it was proved in [9] that $f \in \mathcal{L}_1$ implies that

$$(1.12) \qquad \|f - f_N\| \lesssim N^{-1/6}.$$

This rate was improved to $N^{-11/62}$ in [15], but, on the other hand, it was shown in [19] that for a particular dictionary, there exists $f \in \mathcal{L}_1$ such that

$$(1.13) \qquad \|f - f_N\| \gtrsim N^{-0.27}.$$

When compared with (1.11), we see that the PGA is far from being optimal.

The RGA, OGA and SPA behave somewhat better: it was proven in [14] for the RGA and SPA and in [9] for the OGA that one has

$$(1.14) \qquad \|f - f_N\| \lesssim \|f\|_{\mathcal{L}_1} N^{-1/2}$$

for all $f \in \mathcal{L}_1$.

For each of these algorithms, it is known that the convergence rate $N^{-1/2}$ cannot generally be improved, even for functions which admit a very sparse expansion in the dictionary $\mathcal{D}$ (see [9] for such a result with a function being the sum of two elements of $\mathcal{D}$).

At this point, some remarks are in order regarding the meaning of the condition $f \in \mathcal{L}_1$ for some concrete dictionaries. A commonly made statement is that greedy algorithms break the *curse of dimensionality*, in that the rate $N^{-1/2}$ is independent of the dimension $d$ of the variable space for $f$ and only relies on the assumption that $f \in \mathcal{L}_1$. This is not exactly true since, in practice, the condition that $f \in \mathcal{L}_1$ becomes more and more stringent as $d$ grows. For instance, in the case where we work in the Hilbert space $\mathcal{H} := L_2([0,1]^d)$ and where $\mathcal{D}$ is a *wavelet basis* $(\psi_\lambda)$, it is known that the smoothness property which ensures that $f \in \mathcal{L}_1$ is that $f$ should belong to the Besov space $B_1^s(L_1)$ with $s = d/2$, which roughly means that $f$ has all of its derivatives of order less than or equal to $d/2$ in $L_1$ (see [8] for the characterization of Besov spaces by the properties of wavelet coefficients). Moreover, for this to hold, it is required that the dual wavelets $\tilde{\psi}_\lambda$ have at least $d/2 - 1$ vanishing moments. Another instance is the case where $\mathcal{D}$ consists of sigmoidal functions of the type $\sigma(v \cdot x - w)$, where $\sigma$ is a fixed univariate function, $v$ is an arbitrary vector in $\mathbb{R}^d$, and $w$ is an arbitrary real number. For such dictionaries, it was proved in [4] that a sufficient condition to have $f \in \mathcal{L}_1$ is the convergence of $\int |\omega| |\mathcal{F}f(\omega)| d\omega$, where $\mathcal{F}$ is the Fourier operator. This integrability condition requires a larger amount of decay on the Fourier transform $\mathcal{F}f$ as $d$ grows. Assuming that $f \in \mathcal{L}_1$ is therefore more and more restrictive as $d$ grows. Similar remarks also hold for other dictionaries (hyperbolic wavelets, Gabor functions, etc.).

1.2. *Results of this paper.* The discussion of the previous section points to a significant weakness in the present theory of greedy algorithms, in that there are no viable bounds for the performance of greedy algorithms for general functions $f \in \mathcal{H}$. This is a severe impediment in some application domains (such as learning theory) where there is no a priori knowledge indicating that the target function is in $\mathcal{L}_1$. One of the main contributions of the present paper is to provide error bounds for the performance of greedy algorithms for general functions $f \in \mathcal{H}$. We shall focus our attention on the OGA and RGA, which, as explained above, have better convergence properties in $\mathcal{L}_1$ than the PGA. We shall consider the specific version of the RGA in which $\alpha_k$ is fixed at $(1 - \lambda/k)_+$, for some fixed $\lambda \geq 1$, and then $(\beta_k, g_k)$ are optimized.

Inspection of the proofs in our paper shows that all further approximation results proved for this version of the RGA also hold for *any* greedy algorithm such that

$$(1.15) \qquad \|f - f_k\| \leq \min_{\beta, g} \|f - \alpha_k f_{k-1} + \beta g\|,$$

irrespective of how $f_k$ is defined. In particular, they hold for the more general version of the RGA defined by (1.8), as well as for the SPA.

In Section 2, we introduce both algorithms and recall the optimal approximation rate $N^{-1/2}$ when the target function $f$ is in $\mathcal{L}_1$. Later in this section, we develop a technique based on the interpolation of operators that provides convergence rates $N^{-s}$, $0 < s < 1/2$, whenever $f$ belongs to a certain intermediate space between $\mathcal{L}_1$ and the Hilbert space $\mathcal{H}$. Namely, we use the spaces

$$(1.16) \qquad \mathcal{B}_p := [\mathcal{H}, \mathcal{L}_1]_{\theta, \infty}, \qquad \theta := 2/p - 1, \ 1 < p < 2,$$

which are the real interpolation spaces between $\mathcal{H}$ and $\mathcal{L}_1$. We show that if $f \in \mathcal{B}_p$, then the OGA and RGA, when applied to $f$, provide approximation rates $CN^{-s}$ with $s := \theta/2 = 1/p - 1/2$. Thus, if we set $\mathcal{B}_1 = \mathcal{L}_1$, then these spaces provide a full range of approximation rates for greedy algorithms. Recall, as discussed previously, for general dictionaries, greedy algorithms will not provide convergence rates better than $N^{-1/2}$ for even the simplest of functions. The results we obtain are optimal in the sense that we recover the best possible convergence rate in the case where the dictionary is an orthonormal basis. For an arbitrary target function $f \in \mathcal{H}$, convergence of the OGA and RGA holds without rate. Finally, we conclude the section by discussing several issues related to the numerical implementation of these greedy algorithms. In particular, we consider the effect of reducing the dictionary $\mathcal{D}$ to a finite sub-dictionary.

In Section 3, we consider the learning problem under the assumption that the data $\mathbf{y} := (y_1, \ldots, y_n)$ are bounded in absolute value by some fixed constant $B$. Our estimator is built on the application of the OGA or RGA to the noisy data $\mathbf{y}$ in the Hilbert space defined by the empirical norm

$$\|f\|_n := \frac{1}{n} \sum_{i=1}^{n} |f(x_i)|^2 \tag{1.17}$$

and its associated inner product. At each step $k$, the algorithm generates an approximation $\hat{f}_k$ to the data. Our estimator is defined by

$$\hat{f} := T\hat{f}_{k^*}, \tag{1.18}$$

where

$$Tx := T_B x := \min\{B, |x|\}\operatorname{sgn}(x) \tag{1.19}$$

is the truncation operator at level $B$ and the value of $k^*$ is selected by a complexity regularization procedure. The main result for this estimator is (roughly) that when the regression function $f_\rho$ is in $\mathcal{B}_p$ [where this space is defined with respect to the norm $\|u\|^2 := E(|u(x)|^2)$], the estimator has convergence rate

$$E(\|\hat{f} - f_\rho\|^2) \lesssim \left(\frac{n}{\log n}\right)^{-2s/(1+2s)}, \tag{1.20}$$

again with $s := 1/p - 1/2$. In the case where $f_\rho \in \mathcal{L}_1$, we obtain the same result with $p = 1$ and $s = 1/2$. We also show that this estimator is universally consistent.

In order to place these results within the current state of the art of statistical learning theory, let us first remark that similar convergence rates for the denoising and the learning problem could be obtained by a more "brute force" approach involving the selection of a proper subset of $\mathcal{D}$ by complexity regularization with techniques such as those in [2] or in Chapter 12 of [11]. Following, for instance, the general approach of [11], this would typically first require restricting the size of the dictionary $\mathcal{D}$ [usually to be of size $O(n^a)$ for some $a > 1$] and then considering all possible subsets $\Lambda \subset \mathcal{D}$ and spaces $\mathcal{G}_\Lambda := \operatorname{Span}\{g \in \Lambda\}$, each of them defining an estimator

$$\hat{f}_\Lambda := T_B\left(\operatorname*{Arg\,min}_{f \in \mathcal{G}_\Lambda} \|y - f\|_n^2\right). \tag{1.21}$$

The estimator $\hat f$ is then defined as the $\hat f_\Lambda$ which minimizes

(1.22) $$\min_{\Lambda \subset \mathcal{D}} \{\|y - \hat f_\Lambda\|_n^2 + \operatorname{Pen}(\Lambda, n)\},$$

with $\operatorname{Pen}(\Lambda, n)$ a complexity penalty term. The penalty term usually restricts the size of $\Lambda$ to be at most $\mathcal{O}(n)$, but even then, the search is over $O(n^{an})$ subsets. In some other approaches, the sets $\mathcal{G}_\Lambda$ might also be discretized, transforming the subproblem of selecting $\hat f_\Lambda$ into a discrete optimization problem.

The main advantage of using the greedy algorithm in place of (1.22) for constructing the estimator is a dramatic reduction of the computational cost. Indeed, instead of considering all possible subsets $\Lambda \subset \mathcal{D}$, the algorithm only considers the sets $\Lambda_k := \{g_1, \ldots, g_k\}$, $k = 1, \ldots, n$, generated by the empirical greedy algorithm.

This approach was proposed and analyzed in [18] using a version of the RGA in which

(1.23) $$\alpha_k + \beta_k = 1,$$

which implies that the approximation $f_k$ at each iteration stays in the convex hull $\mathcal{C}_1$ of $\mathcal{D}$. The authors established that if $f$ does not belong to $\mathcal{C}_1$, then the RGA converges to its projection onto $\mathcal{C}_1$. In turn, the estimator was proven to converge in the sense of (1.20) to $f_\rho$, with rate $(n/\log n)^{-1/2}$, if $f_\rho$ lies in $\mathcal{C}_1$ and otherwise to its projection onto $\mathcal{C}_1$. In that sense, this procedure is not universally consistent.

One of the main contributions of the present paper is to remove requirements of the type $f_\rho \in \mathcal{L}_1$ when obtaining convergence rates. In the learning context, there is indeed typically no advanced information that would guarantee such restrictions on $f_\rho$. The estimators that we construct for learning are now universally consistent and have provable convergence rates for more general regression functions described by means of interpolation spaces. One of the main ingredients in our analysis of the performance of our greedy algorithms in learning is a powerful exponential concentration inequality which was introduced in [18]. Let us mention that a closely related analysis, which does not, however, involve interpolation spaces, is developed in [5, 13].

The most studied dictionaries in learning theory are in the context of neural networks. In Section 4, we interpret our results in this setting and, in particular, describe the smoothness conditions on a function $f$ which ensure that it belongs to the spaces $\mathcal{L}_1$ or $\mathcal{B}_p$.

Let us finally mention that there exist some natural connections between the greedy algorithms which are discussed in this paper and other numerical techniques for building a sparse approximation in the dictionary based on the minimization of an $\ell_1$ penalized criterion. In the statistical context, these are the celebrated LASSO [12, 23] and LARS [10] algorithms. The relation between $\ell_1$ minimization and greedy selection is particularly transparent in the context of deterministic approximation of a function $f$ in an orthonormal basis: if we consider the problem of minimizing

$$\left\| f - \sum_{g \in \mathcal{D}} d_g g \right\|^2 + t \sum_{g \in \mathcal{D}} |d_g| \tag{1.24}$$

over all choices of sequences $(d_g)$, we see that it amounts to minimizing $|c_g - d_g|^2 + t|d_g|$ for each individual $g$, where $c_g := \langle f, g \rangle$. The solution to this problem is given by the *soft thresholding operator*

$$d_g := \begin{cases} c_g - \dfrac{t}{2} \mathrm{sign}(c_g), & \text{if } |c_g| > \dfrac{t}{2}, \\ 0, & \text{otherwise} \end{cases} \tag{1.25}$$

and is therefore very similar to the results of selecting the largest coefficients of $f$.

**2. Approximation properties.** Let $\mathcal{D}$ be a dictionary in some Hilbert space $\mathcal{H}$, with $\|g\| = 1$ for all $g \in \mathcal{D}$. We recall that, for a given $f \in \mathcal{H}$, the OGA builds embedded approximation spaces

$$V_k := \mathrm{Span}\{g_1, \ldots, g_k\}, \qquad k = 1, 2, \ldots, \tag{2.1}$$

and approximations

$$f_k := P_k f, \tag{2.2}$$

where $P_k$ is the orthogonal projection onto $V_k$. The rule for generating the $g_k$ is as follows. We set $V_0 := \{0\}$, $f_0 := 0$ and $r_0 := f$ and, given $V_{k-1}$, $f_{k-1} = P_{k-1} f$ and $r_{k-1} := f - f_{k-1}$, we define $g_k$ by

$$g_k := \underset{g \in G}{\mathrm{Arg\,max}}\, |\langle r_{k-1}, g \rangle|, \tag{2.3}$$

which defines the new $V_k$, $f_k$ and $r_k$.

In its most general form, the RGA sets

$$f_k = \alpha_k f_{k-1} + \beta_k g_k, \tag{2.4}$$

where $(\alpha_k, \beta_k, g_k)$ are defined according to (1.8). We shall consider a simpler version in which the first parameter is a fixed sequence. Two choices will be considered, namely

$$\alpha_k = 1 - \frac{1}{k} \tag{2.5}$$

and

$$\alpha_k = 1 - \frac{2}{k} \qquad \text{if } k > 1,\ \alpha_1 = 0. \tag{2.6}$$

The two other parameters are optimized according to

$$(2.7) \qquad (\beta_k, g_k) := \underset{(\beta,g) \in \mathbb{R} \times \mathcal{D}}{\operatorname{Arg\,min}} \|f - \alpha_k f_{k-1} - \beta g\|.$$

Since

$$(2.8) \quad \|f - \alpha_k f_{k-1} - \beta g\|^2 = \beta^2 - 2\beta \langle f - \alpha_k f_{k-1}, g \rangle + \|f - \alpha_k f_{k-1}\|^2,$$

it is readily seen that $(\beta_k, g_k)$ are given explicitly by

$$(2.9) \qquad \beta_k = \langle f - \alpha_k f_{k-1}, g_k \rangle$$

and

$$(2.10) \qquad g_k := \underset{g \in \mathcal{D}}{\operatorname{Arg\,max}} |\langle f - \alpha_k f_{k-1}, g \rangle|.$$

Therefore, from a computational point of view, this RGA is very similar to the PGA.

We denote by $\mathcal{L}_p$ the functions $f$ which admit a converging expansion $f = \sum c_g g$ with $\sum |c_g|^p < +\infty$ and we write $\|f\|_{\mathcal{L}_p} = \inf \|(c_g)\|_{\ell_p}$, where the infimum is taken over all such expansions. In a similar way, we consider the spaces $w\mathcal{L}_p$ corresponding to expansions which are in the weak space $w\ell_p$. We denote by $\sigma_N(f)$ the best $N$-term approximation error in the $\mathcal{H}$ norm for $f$ and for any $s > 0$, we define the approximation space

$$(2.11) \qquad \mathcal{A}^s := \{f \in \mathcal{H} : \sigma_N(f) \leq MN^{-s}, N = 1, 2, \ldots\}.$$

Finally, we denote by $\mathcal{G}^s$ the set of functions $f$ such that the greedy algorithm under consideration converges with rate $\|r_N\| \lesssim N^{-s}$, so that, obviously, $\mathcal{G}^s \subset \mathcal{A}^s$.

In the case where $\mathcal{D}$ is an orthonormal basis, the space $\mathcal{A}^s$ contains the space $\mathcal{L}_p$ with $1/p = 1/2 + s$ and, in fact, actually coincides with the weak versions $w\mathcal{L}_p$ of these spaces. In those cases, an algorithm for building a best (or near best) $N$-term approximation is simply to keep the $N$ largest coefficients of $f$ and discard the others. The best $N$-term approximation is also obtained by the orthogonal greedy algorithm so that, obviously, $\mathcal{A}^s = \mathcal{G}^s$.

2.1. *Approximation of $\mathcal{L}_1$ functions.* In this section, we recall, for convenience, the approximation properties of the OGA and RGA for functions $f \in \mathcal{L}_1$. We first recall the result obtained in [9] for the OGA. We shall make use of the following fact: if $f, g \in \mathcal{H}$ with $\|g\| = 1$, then $\langle f, g \rangle g$ is the best approximation to $f$ from the one-dimensional space generated by $g$ and

$$(2.12) \qquad \|f - \langle f, g \rangle g\|^2 = \|f\|^2 - |\langle f, g \rangle|^2.$$

THEOREM 2.1. *For all $f \in \mathcal{L}_1$, the error of the OGA satisfies*

(2.13) $$\|r_N\| \leq \|f\|_{\mathcal{L}_1}(N+1)^{-1/2}, \qquad N = 1, 2, \ldots.$$

PROOF. Since $f_k$ is the best approximation to $f$ from $V_k$, we have, from (2.12),

(2.14) $$\|r_k\|^2 \leq \|r_{k-1} - \langle r_{k-1}, g_k\rangle g_k\|^2 = \|r_{k-1}\|^2 - |\langle r_{k-1}, g_k\rangle|^2,$$

with equality in the case where $g_k$ is orthogonal to $V_{k-1}$. Since $r_{k-1}$ is orthogonal to $f_{k-1}$, we have

(2.15) $$\|r_{k-1}\|^2 = \langle r_{k-1}, f\rangle \leq \|f\|_{\mathcal{L}_1} \sup_{g \in \mathcal{D}} |\langle r_{k-1}, g\rangle| = \|f\|_{\mathcal{L}_1} |\langle r_{k-1}, g_k\rangle|,$$

which, combined with (2.14), gives the reduction property

(2.16) $$\|r_k\|^2 \leq \|r_{k-1}\|^2 (1 - \|r_{k-1}\|^2 \|f\|_{\mathcal{L}_1}^{-2}).$$

We also know that $\|r_1\| \leq \|r_0\| = \|f\| \leq \|f\|_{\mathcal{L}_1}$.

We then check by induction that a decreasing sequence $(a_n)_{n \geq 0}$ of nonnegative numbers which satisfy $a_0 \leq M$ and $a_k \leq a_{k-1}(1 - \frac{a_{k-1}}{M})$ for all $k > 0$ has the decay property $a_n \leq \frac{M}{n+1}$ for all $n \geq 0$. Indeed, assuming that $a_{n-1} \leq \frac{M}{n}$ for some $n > 0$, then either $a_{n-1} \leq \frac{M}{n+1}$, so that $a_n \leq \frac{M}{n+1}$, or else $a_{n-1} \geq \frac{M}{n+1}$, so that

(2.17) $$a_n \leq \frac{M}{n}\left(1 - \frac{1}{n+1}\right) = \frac{M}{n+1}.$$

The result follows by applying this to $a_k = \|r_k\|^2$ and $M := \|f\|_{\mathcal{L}_1}^2$, since we indeed have

(2.18) $$a_0 = \|f\|^2 \leq \|f\|_{\mathcal{L}_1}^2. \qquad \square$$

We now turn to the RGA, for which we shall prove a slightly stronger property.

THEOREM 2.2. *For all $f \in \mathcal{L}_1$, the error of the RGA using* (2.5) *satisfies*

(2.19) $$\|r_N\| \leq (\|f\|_{\mathcal{L}_1}^2 - \|f\|^2)^{1/2} N^{-1/2}, \qquad N = 1, 2, \ldots.$$

PROOF. From the definition of the RGA, we see that the sequence $f_k$ remains unchanged if the dictionary $\mathcal{D}$ is symmetrized by including the sequence $(-g)_{g \in \mathcal{D}}$. Under this assumption, since $f \in \mathcal{L}_1$, for any $\varepsilon > 0$, we can expand $f$ according to

(2.20) $$f = \sum_{g \in \mathcal{D}} b_g g,$$

where all of the $b_g$ are nonnegative and satisfy

(2.21) $$\sum_{g \in \mathcal{D}} b_g = \|f\|_{\mathcal{L}_1} + \delta,$$

with $0 \leq \delta \leq \varepsilon$. According to (2.7), we have, for all $\beta \in \mathbb{R}$ and all $g \in \mathcal{D}$,

$$\|r_k\|^2 \leq \|f - \alpha_k f_{k-1} - \beta g\|^2$$
$$= \left\| \alpha_k r_{k-1} + \frac{1}{k} f - \beta g \right\|^2$$
$$= \alpha_k^2 \|r_{k-1}\|^2 - 2\alpha_k \left\langle r_{k-1}, \frac{1}{k} f - \beta g \right\rangle + \left\| \frac{1}{k} f - \beta g \right\|^2$$
$$= \alpha_k^2 \|r_{k-1}\|^2 - 2\alpha_k \left\langle r_{k-1}, \frac{1}{k} f - \beta g \right\rangle + \frac{1}{k^2} \|f\|^2 + \beta^2 - \frac{2\beta}{k} \langle f, g \rangle.$$

This inequality holds for all $g \in \mathcal{D}$, so it also holds on the average with weights $\frac{b_g}{\sum_{g \in \mathcal{D}} b_g}$, which gives, for the particular value $\beta = \frac{1}{k}(\|f\|_{\mathcal{L}_1} + \delta)$,

(2.22) $$\|r_k\|^2 \leq \alpha_k^2 \|r_{k-1}\|^2 - \frac{1}{k^2} \|f\|^2 + \beta^2.$$

Therefore, letting $\varepsilon$ tend to 0, we obtain

(2.23) $$\|r_k\|^2 \leq \left(1 - \frac{1}{k}\right)^2 \|r_{k-1}\|^2 + \frac{1}{k^2}(\|f\|_{\mathcal{L}_1}^2 - \|f\|^2).$$

We now check by induction that a sequence $(a_k)_{k>0}$ of positive numbers such that $a_1 \leq M$ and $a_k \leq (1 - \frac{1}{k})^2 a_{k-1} + \frac{1}{k^2} M$ for all $k > 0$ has the decay property $a_n \leq \frac{M}{n}$ for all $n > 0$. Indeed, assuming that $a_{n-1} \leq \frac{M}{n-1}$, we write

$$a_n - \frac{M}{n} \leq \left(1 - \frac{1}{n}\right)^2 \frac{M}{n-1} + \frac{1}{n^2} M - \frac{M}{n}$$
$$= M\left(\frac{n-1}{n^2} + \frac{1}{n^2} - \frac{1}{n}\right) = 0.$$

The result follows by applying this to $a_k = \|r_k\|^2$ and $M := \|f\|_{\mathcal{L}_1}^2 - \|f\|^2$ since (2.23) also implies that $a_1 \leq M$. □

The above results show that for both OGA and RGA, we have

(2.24) $$\mathcal{L}_1 \subset \mathcal{G}^{1/2} \subset \mathcal{A}^{1/2}.$$

From this, it also follows that $w\mathcal{L}_p \subset \mathcal{A}^s$ with $s = 1/p - 1/2$ when $p < 1$. Indeed, from the definition of $w\mathcal{L}_p$, any function $f$ in this space can be written as $f = \sum_{j=1}^{\infty} c_j g_j$, with each $g_j \in \mathcal{D}$ and the coefficients $c_j$ decreasing in absolute value and satisfying $|c_j| \leq M j^{-1/p}$, $j \geq 1$, with $M := \|f\|_{w\mathcal{L}_p}$. Therefore, $f = f_a + f_b$ with $f_a := \sum_{j=1}^{N} c_j g_j$ and $\|f_b\|_{\mathcal{L}_1} \leq C_p M N^{1-1/p}$. It follows

from Theorem 2.1, for example, that $f_b$ can be approximated by an $N$-term expansion obtained by the greedy algorithm with accuracy $\|f_b - P_N f_b\| \leq C_p M N^{1/2-1/p} = N^{-s}$ and therefore, by taking $f_a + P_N f_b$ as a $2N$-term approximant of $f$, we obtain that $f \in \mathcal{A}^s$. Observe, however, that this does not mean that $f \in \mathcal{G}^s$ in the sense that we have not proven that the greedy algorithm converges with rate $N^{-s}$ when applied to $f$. It is actually shown in [9] that there exist simple dictionaries such that the greedy algorithm does not converge faster than $N^{-1/2}$, even for functions $f$ which are in $w\mathcal{L}_p$ for all $p > 0$.

2.2. *Approximation of general functions.* We now want to study the behavior of the OGA and RGA when the function $f \in \mathcal{H}$ is more general, in the sense that it is less sparse than being in $\mathcal{L}_1$. The simplest way of expressing this would seem to be by considering the spaces $\mathcal{L}_p$ or $w\mathcal{L}_p$ with $1 < p < 2$. However, for general dictionaries, these spaces are not well adapted, since $\|f\|_{\mathcal{L}_p}$ does not control the Hilbert norm $\|f\|$.

Instead, we shall consider the real interpolation spaces

$$\mathcal{B}_p = [\mathcal{H}, \mathcal{L}_1]_{\theta,\infty}, \qquad 0 < \theta < 1, \tag{2.25}$$

again with $p$ defined by $1/p = \theta + (1-\theta)/2 = (1+\theta)/2$. Recall that $f \in [X,Y]_{\theta,\infty}$ if and only if for all $t > 0$, we have

$$K(f,t) \leq Ct^\theta, \tag{2.26}$$

where

$$K(f,t) := K(f,t,X,Y) := \inf_{h \in Y}\{\|f-h\|_X + t\|h\|_Y\} \tag{2.27}$$

is the so-called *K-functional*. In other words, $f$ can be decomposed into $f = f_X + f_Y$, with

$$\|f_X\|_X + t\|f_Y\|_Y \leq Ct^\theta. \tag{2.28}$$

The smallest $C$ such that the above holds defines a norm for $Z = [X,Y]_{\theta,\infty}$. We refer to [7] or [6] for an introduction to interpolation spaces. The space $\mathcal{B}_p$ coincides with $w\mathcal{L}_p$ in the case where $\mathcal{D}$ is an orthonormal system, but may differ from it for a more general dictionary.

The main result of this section is that for both the OGA and the RGA,

$$\|r_N\| \leq C_0 K(f, N^{-1/2}, \mathcal{H}, \mathcal{L}_1), \qquad N = 1, 2, \ldots, \tag{2.29}$$

so that, according to (2.26), $f \in \mathcal{B}_p$ implies the rate of decay $\|r_N\| \lesssim N^{-\theta/2}$. Note that if $f_N$ were obtained as the action on $f$ of a continuous linear operator $L_N$ from $\mathcal{H}$ onto itself such that $\|L_N\| \leq C$ with $C$ independent of $N$, then we could write, for any $h \in \mathcal{L}_1$,

$$\begin{aligned}\|f - f_N\| &\leq \|(I - L_N)[f-h]\| + \|h - L_N h\| \\ &\lesssim \|f - h\| + \|h\|_{\mathcal{L}_1} N^{-1/2},\end{aligned} \tag{2.30}$$

so that (2.29) would follow by minimizing over $h \in \mathcal{L}_1$. However, $f_N$ is obtained by a highly nonlinear algorithm and it is therefore quite remarkable that (2.29) still holds. We first prove this for the OGA.

THEOREM 2.3. *For all $f \in \mathcal{H}$ and any $h \in \mathcal{L}_1$, the error of the OGA satisfies*

$$\|r_N\|^2 \leq \|f - h\|^2 + 4\|h\|_{\mathcal{L}_1}^2 N^{-1}, \qquad N = 1, 2, \ldots, \tag{2.31}$$

*and therefore*

$$\begin{aligned}\|r_N\| &\leq K(f, 2N^{-1/2}, \mathcal{H}, \mathcal{L}_1) \\ &\leq 2K(f, N^{-1/2}, \mathcal{H}, \mathcal{L}_1), \qquad N = 1, 2, \ldots.\end{aligned} \tag{2.32}$$

PROOF. Fix an arbitrary $f \in \mathcal{H}$. For any $h \in \mathcal{L}_1$, we write

$$\|r_{k-1}\|^2 = \langle r_{k-1}, h + f - h \rangle \leq \|h\|_{\mathcal{L}_1} |\langle r_{k-1}, g_k \rangle| + \|r_{k-1}\| \|f - h\|, \tag{2.33}$$

from which it follows that

$$\|r_{k-1}\|^2 \leq \|h\|_{\mathcal{L}_1} |\langle r_{k-1}, g_k \rangle| + \tfrac{1}{2}(\|r_{k-1}\|^2 + \|f - h\|^2). \tag{2.34}$$

Therefore, using the shorthand notation $a_k := \|r_k\|^2 - \|f - h\|^2$, we have

$$|\langle r_{k-1}, g_k \rangle| \geq \frac{a_{k-1}}{2\|h\|_{\mathcal{L}_1}}. \tag{2.35}$$

Note that if for some $k_0$, we have $\|r_{k_0-1}\| \leq \|f - h\|$, then the theorem holds trivially for all $N \geq k_0 - 1$. We therefore assume that $a_{k-1}$ is positive, so that we can write

$$|\langle r_{k-1}, g_k \rangle|^2 \geq \frac{a_{k-1}^2}{4\|h\|_{\mathcal{L}_1}^2}. \tag{2.36}$$

From (2.14), we therefore obtain

$$\|r_k\|^2 \leq \|r_{k-1}\|^2 - \frac{a_{k-1}^2}{4\|h\|_{\mathcal{L}_1}^2}, \tag{2.37}$$

which, by subtracting $\|f - h\|^2$, gives

$$a_k \leq a_{k-1}\left(1 - \frac{a_{k-1}}{4\|h\|_{\mathcal{L}_1}^2}\right). \tag{2.38}$$

As in the proof of Theorem 2.1, we can conclude that $a_N \leq 4\|h\|_{\mathcal{L}_1}^2 N^{-1}$, provided that we initially have $a_1 \leq 4\|h\|_{\mathcal{L}_1}^2$. In order to check this initial condition, we remark that either $a_0 \leq 4\|h\|_{\mathcal{L}_1}^2$, so that the same holds for $a_1$, or $a_0 \geq 4\|h\|_{\mathcal{L}_1}^2$, in which case $a_1 \leq 0$ according to (2.38), which means that we are already in the trivial case $\|r_1\| \leq \|f - h\|$ for which there is nothing to prove. We have therefore obtained (2.31) and (2.32) follows by taking the square root. $\square$

We next treat the case of the RGA, for which we have a slightly different result. In this result, we use the second choice (2.6) for the sequence $\alpha_k$ in order to obtain a multiplicative constant equal to 1 in the term $\|f - h\|^2$ appearing on the right side of the quadratic bound. This will be important in the learning application. We also give the nonquadratic bound with the first choice (2.5) since it gives a slightly better result than by taking the square root of the quadratic bound based on (2.6).

THEOREM 2.4. *For all $f \in \mathcal{H}$ and any $h \in \mathcal{L}_1$, the error of the RGA using* (2.6) *satisfies*

$$(2.39) \quad \|r_N\|^2 \leq \|f - h\|^2 + 4(\|h\|_{\mathcal{L}_1}^2 - \|h\|^2)N^{-1}, \qquad N = 1, 2, \ldots,$$

*and therefore*

$$(2.40) \quad \begin{aligned} \|r_N\| &\leq K(f, 2N^{-1/2}, \mathcal{H}, \mathcal{L}_1) \\ &\leq 2K(f, N^{-1/2}, \mathcal{H}, \mathcal{L}_1), \qquad N = 1, 2, \ldots. \end{aligned}$$

*Using the first choice* (2.5), *the error satisfies*

$$(2.41) \quad \|r_N\| \leq \|f - h\| + (\|h\|_{\mathcal{L}_1}^2 - \|h\|^2)^{1/2} N^{-1/2}, \qquad N = 1, 2, \ldots,$$

*and therefore*

$$(2.42) \quad \|r_N\| \leq K(f, N^{-1/2}, \mathcal{H}, \mathcal{L}_1), \qquad N = 1, 2, \ldots.$$

PROOF. Fix $f \in \mathcal{H}$ and let $h \in \mathcal{L}_1$ be arbitrary. Similarly to the proof of Theorem 2.2, for any $\varepsilon > 0$, we can expand $h$ as

$$(2.43) \quad h = \sum_{g \in \mathcal{D}} b_g g,$$

where all of the $b_g$ are nonnegative and satisfy

$$(2.44) \quad \sum_{g \in \mathcal{D}} b_g = \|h\|_{\mathcal{L}_1} + \delta,$$

with $0 \leq \delta \leq \varepsilon$. Using the notation $\bar{\alpha}_k = 1 - \alpha_k$, we have, for all $\beta \in \mathbb{R}$ and all $g \in \mathcal{D}$,

$$\begin{aligned} \|r_k\|^2 &\leq \|f - \alpha_k f_{k-1} - \beta g\|^2 \\ &= \|\alpha_k r_{k-1} + \bar{\alpha}_k f - \beta g\|^2 \\ &= \alpha_k^2 \|r_{k-1}\|^2 - 2\alpha_k \langle r_{k-1}, \bar{\alpha}_k f - \beta g \rangle + \|\bar{\alpha}_k f - \beta g\|^2 \end{aligned}$$

$$= \alpha_k^2 \|r_{k-1}\|^2 - 2\alpha_k \langle r_{k-1}, \bar{\alpha}_k f - \beta g\rangle + \|\bar{\alpha}_k(f-h) + \bar{\alpha}_k h - \beta g\|^2$$
$$= \alpha_k^2 \|r_{k-1}\|^2 - 2\alpha_k \langle r_{k-1}, \bar{\alpha}_k f - \beta g\rangle + \bar{\alpha}_k^2 \|f-h\|^2$$
$$+ 2\bar{\alpha}_k \langle f-h, \bar{\alpha}_k h - \beta g\rangle + \bar{\alpha}_k^2 \|h\|^2 - 2\beta \bar{\alpha}_k \langle h, g\rangle + \beta^2.$$

This inequality holds for all $g \in \mathcal{D}$, so it also holds on the average with weights $\frac{b_g}{\sum_{g\in\mathcal{D}} b_g}$, which gives, for the particular value $\beta = \bar{\alpha}_k(\|h\|_{\mathcal{L}_1} + \delta)$,

$$\|r_k\|^2 \le \alpha_k^2 \|r_{k-1}\|^2 - 2\alpha_k \bar{\alpha}_k \langle r_{k-1}, f-h\rangle + \bar{\alpha}_k^2 \|f-h\|^2 - \bar{\alpha}_k^2 \|h\|^2 + \beta^2$$
$$= \|\alpha_k r_{k-1} - \bar{\alpha}_k(f-h)\|^2 - \bar{\alpha}_k^2 \|h\|^2 + \beta^2.$$

Letting $\varepsilon$ tend to 0 and using the notation $M := \|h\|_{\mathcal{L}_1}^2 - \|h\|^2$, we thus obtain

$$(2.45) \qquad \|r_k\|^2 \le (\alpha_k \|r_{k-1}\| + \bar{\alpha}_k \|f-h\|)^2 + \bar{\alpha}_k^2 M.$$

Note that this immediately implies the validity of (2.39) and (2.41) at $N=1$, using the fact that $\alpha_1 = 0$ for both choices (2.5) and (2.6). We next proceed by induction, assuming that these bounds hold at $k-1$.

For the proof of (2.41), we derive from (2.45) that

$$\|r_k\|^2 \le \left(\alpha_k\left(\|f-h\| + \left(\frac{M}{k-1}\right)^{1/2}\right) + \bar{\alpha}_k \|f-h\|\right)^2 + \bar{\alpha}_k^2 M$$
$$= \left(\|f-h\| + \alpha_k\left(\frac{M}{k-1}\right)^{1/2}\right)^2 + \bar{\alpha}_k^2 M$$
$$= \|f-h\|^2 + 2M^{1/2}\|f-h\|\frac{1-1/k}{\sqrt{k-1}} + M\left(\frac{(1-1/k)^2}{k-1} + \frac{1}{k^2}\right)$$
$$= \|f-h\|^2 + 2M^{1/2}\|f-h\|\frac{\sqrt{k-1}}{k} + \frac{M}{k}$$
$$\le \|f-h\|^2 + 2\left(\frac{M}{k}\right)^{1/2}\|f-h\| + \frac{M}{k}$$
$$= \left(\|f-h\| + \left(\frac{M}{k}\right)^{1/2}\right)^2,$$

which is the desired bound at $k$.

For the proof of (2.39), we derive from (2.45) that

$$(2.46) \qquad \begin{aligned}\|r_k\|^2 &\le \alpha_k^2 \|r_{k-1}\|^2 + 2\alpha_k \bar{\alpha}_k \|r_{k-1}\|\|f-h\| \\ &\quad + \bar{\alpha}_k^2 \|f-h\|^2 + \bar{\alpha}_k^2 M.\end{aligned}$$

Noting that $2\alpha_k \bar{\alpha}_k \|r_{k-1}\|\|f-h\| \le \alpha_k \bar{\alpha}_k(\|r_{k-1}\|^2 + \|f-h\|^2)$, we obtain

$$(2.47) \qquad \|r_k\|^2 \le \alpha_k \|r_{k-1}\|^2 + \bar{\alpha}_k \|f-h\|^2 + \bar{\alpha}_k^2 M$$

and therefore

$$\|r_k\|^2 - \|f - h\|^2 \tag{2.48}$$
$$\leq \alpha_k(\|r_{k-1}\|^2 - \|f - h\|^2) + \bar{\alpha}_k^2 M.$$

Assuming that (2.39) holds at $N = k - 1$, we thus obtain

$$\|r_k\|^2 - \|f - h\|^2 \leq \left(1 - \frac{2}{k}\right)\frac{4M}{k-1} + \frac{4M}{k^2} \leq \frac{4M}{k}. \tag{2.49}$$

We have thus obtained (2.39), and (2.40) follows by taking the square root. □

REMARK 2.5. More generally, we can consider the RGA with the choice

$$\alpha_k = \left(1 - \frac{\lambda}{k}\right)_+, \tag{2.50}$$

where $\lambda$ is some fixed parameter. If $\lambda > 1$, an argument similar that used in the proof of (2.39) gives the general estimate

$$\|r_N\|^2 \leq \|f - h\|^2 + C(\|h\|_{\mathcal{L}_1}^2 - \|h\|^2)N^{-1}, \tag{2.51}$$

with $C = \frac{\lambda^2}{\lambda - 1}$. A specific feature of this estimate is that the term $\|f - h\|^2$ has a multiplicative constant equal to 1, which will be of critical importance in the learning context. We do not know if such an estimate can be obtained if $0 < \lambda \leq 1$.

An immediate consequence of Theorems 2.3 and 2.4 combined with the definition of the $\mathcal{B}_p$ spaces [see (2.26)] is a rate of convergence of the OGA and RGA for the functions in $\mathcal{B}_p$.

COROLLARY 2.6. *For all $f \in \mathcal{B}_p$, the approximation errors for both the OGA and RGA satisfy the decay bound*

$$\|r_N\| \lesssim \|f\|_{\mathcal{B}_p} N^{-s}, \tag{2.52}$$

*with $s = 1/p - 1/2$. Therefore, we have $\mathcal{B}_p \subset \mathcal{G}^s \subset \mathcal{A}^s$.*

In addition, when $\mathcal{D}$ is a complete family in $\mathcal{H}$, we know that $\mathcal{L}_1$ is dense in $\mathcal{H}$, so

$$\lim_{t \to 0} K(f, t, \mathcal{H}, \mathcal{L}_1) = 0 \tag{2.53}$$

for any $f \in \mathcal{H}$. This implies the following corollary.

COROLLARY 2.7. *For any $f \in \mathcal{H}$, the approximation error $\|r_N\|$ tends to zero as $N \to +\infty$ for both the OGA and RGA.*

2.3. *Greedy algorithms with a truncated dictionary.* In concrete applications, it is not possible to evaluate the supremum of $|\langle r_{k-1}, g \rangle|$ over the whole dictionary $\mathcal{D}$, but only over a finite subset of it. For applications in learning theory, it will also be useful that the size of this subset has at most polynomial growth in the number of samples $n$. We therefore introduce a fixed exhaustion of $\mathcal{D}$,

$$\mathcal{D}_1 \subset \mathcal{D}_2 \subset \cdots \subset \mathcal{D}, \tag{2.54}$$

with $\#(\mathcal{D}_m) = m$. The analysis we present in this section is similar to that in [22]. We are now interested in the functions $f$ which can be approximated to a certain accuracy by application of the OGA only using the elements of $\mathcal{D}_m$. For this purpose, we first introduce the space $\mathcal{L}_1(\mathcal{D}_m)$ of those functions in $\mathrm{Span}(\mathcal{D}_m)$ equipped with the (minimal) $\ell_1$ norm of the coefficients. We next define, for $r > 0$, the space $\mathcal{L}_{1,r}$ to be the set of all functions $f$ such that, for all $m$, there exists $h$ (depending on $m$) such that

$$\|h\|_{\mathcal{L}_1(\mathcal{D}_m)} \leq C \tag{2.55}$$

and

$$\|f - h\| \leq C m^{-r}. \tag{2.56}$$

The smallest constant $C$ such that this holds defines a norm for $\mathcal{L}_{1,r}$. In order to understand how these spaces are related to the space $\mathcal{L}_1$ for the whole dictionary, consider the example where $\mathcal{D}$ is a Schauder basis and consider the decomposition of $f$ into

$$\begin{aligned}f &= \sum_{g \in \mathcal{D}_m} c_g g + \sum_{g \notin \mathcal{D}_m} c_g g \\ &= h + f - h.\end{aligned} \tag{2.57}$$

It is then obvious that $\|h\|_{\mathcal{L}_1(\mathcal{D}_m)} \leq \|f\|_{\mathcal{L}_1}$. Therefore, a sufficient condition for $f$ to be in $\mathcal{L}_{1,r}$ is $f \in \mathcal{L}_1$ and its tail $\|\sum_{g \notin \mathcal{D}_m} c_g g\|$ decays like $m^{-r}$.

Application of Theorems 2.3 and 2.4 shows us that if we apply the OGA or RGA with the restricted dictionary and if the target function $f$ is in $\mathcal{L}_{1,r}$, then we have

$$\|r_k\| \leq C_0 \|f\|_{\mathcal{L}_{1,r}} (k^{-1/2} + m^{-r}), \tag{2.58}$$

where $C_0$ is an absolute constant [$C_0 = 2$ for OGA and $C_0 = 1$ for RGA with choice (2.5)].

In a similar manner, we can introduce the interpolation space

$$\mathcal{B}_{p,r} := [\mathcal{H}, \mathcal{L}_{1,r}]_{\theta, \infty}, \tag{2.59}$$

again with $1/p = (1 + \theta)/2$. From the definition of interpolation spaces, if $f \in \mathcal{B}_{p,r}$, then for all $t > 0$, there exists $\tilde{f} \in \mathcal{L}_{1,r}$ such that

$$\|\tilde{f}\|_{\mathcal{L}_{1,r}} \leq \|f\|_{\mathcal{B}_{p,r}} t^{\theta - 1} \tag{2.60}$$

and

$$\|f - \tilde{f}\| \leq \|f\|_{\mathcal{B}_{p,r}} t^\theta. \tag{2.61}$$

We also know that for all $m$, there exists $h$ (depending on $m$) such that

$$\|h\|_{\mathcal{L}_1(\mathcal{D}_m)} \leq \|\tilde{f}\|_{\mathcal{L}_{1,r}} \leq \|f\|_{\mathcal{B}_{p,r}} t^{\theta-1} \tag{2.62}$$

and

$$\begin{aligned}\|\tilde{f} - h\| &\leq \|\tilde{f}\|_{\mathcal{L}_{1,r}} m^{-r} \\ &\leq \|f\|_{\mathcal{B}_{p,r}} t^{\theta-1} m^{-r},\end{aligned} \tag{2.63}$$

so that, by the triangle inequality,

$$\|f - h\| \leq \|f\|_{\mathcal{B}_{p,r}}(t^\theta + t^{\theta-1} m^{-r}). \tag{2.64}$$

Application of Theorems 2.3 and 2.4 shows us that if we apply the OGA or RGA with the restricted dictionary and if the target function $f$ is in $\mathcal{B}_{p,r}$, then we have, for any $t > 0$,

$$\|r_k\| \leq C_0 \|f\|_{\mathcal{B}_{p,r}}(t^{\theta-1} k^{-1/2} + t^\theta + t^{\theta-1} m^{-r}). \tag{2.65}$$

In particular, taking $t = k^{-1/2}$ and noting that $\theta = 2s$ gives

$$\|r_k\| \leq C_0 \|f\|_{\mathcal{B}_{p,r}}(k^{-s} + k^{1/2-s} m^{-r}). \tag{2.66}$$

We therefore recover the rate of Corollary 2.6 up to an additive perturbation which tends to 0 as $m \to +\infty$.

Let us conclude this section by making some remarks on the spaces $\mathcal{B}_{p,r}$. These spaces should be viewed as being slightly smaller than the spaces $\mathcal{B}_p$. The smaller the value of $r > 0$, the smaller the distinction between $\mathcal{B}_p$ and $\mathcal{B}_{p,r}$. Also, note that the classes $\mathcal{B}_{p,r}$ depend very much on how we exhaust the dictionary $\mathcal{D}$. For example, if $\mathcal{D} = B_0 \cup B_1$ is the union of two bases $B_0$ and $B_1$, then exhausting the elements of $B_0$ faster than those of $B_1$ will result in different classes than if we exhaust those of $B_1$ faster than those of $B_0$. However, in concrete settings, there is usually a natural order in which to exhaust the dictionary.

## 3. Application to statistical learning.

3.1. *Notation and definition of the estimator.* We consider the classical bounded regression problem. We observe $n$ independent realizations $(z_i) = (x_i, y_i)$, $i = 1, \ldots, n$, of an unknown distribution $\rho$ on $Z = X \times Y$. We assume here that the output variable satisfies almost surely

$$|y| \leq B, \tag{3.1}$$

where the bound $B$ is known to us. We denote by

$$f_\rho(x) = E(y|x) \tag{3.2}$$

the regression function which minimizes the quadratic risk

$$R(f) := E(|f(x) - y|^2) \tag{3.3}$$

over all functions $f$. For any $f$, we have

$$R(f) - R(f_\rho) = \|f - f_\rho\|^2, \tag{3.4}$$

where we use the notation

$$\|u\|^2 := E(|u(x)|^2) = \|u\|^2_{L_2(\rho_X)}, \tag{3.5}$$

with $\rho_X$ being the marginal probability measure defined on $X$. We are therefore interested in constructing from our data an estimator $\hat{f}$ such that $\|\hat{f} - f_\rho\|^2$ is small. Since $\hat{f}$ depends on the realization of the training sample $\mathbf{z} := (z_i) \in Z^n$, we shall measure the estimation error by the expectation $E(\|\hat{f} - f_\rho\|^2)$ taken with respect to $\rho^n$.

Given our training sample $\mathbf{z}$, we define the empirical norm

$$\|f\|_n^2 := \frac{1}{n} \sum_{i=1}^n |f(x_i)|^2. \tag{3.6}$$

Note that $\|\cdot\|_n$ is the $L_2$ norm with respect to the discrete measure $\nu_\mathbf{x} := \frac{1}{n} \sum_{i=1}^n \delta_{x_i}$, with $\delta_u$ the Dirac measure at $u$. As such, the norm depends on $\mathbf{x} := (x_1, \ldots, x_n)$ and not just on $n$, but we adopt the notation (3.6) to conform with other major works in learning. We view the vector $\mathbf{y} := (y_1, \ldots, y_n)$ as a function $y$ defined on the design $\mathbf{x} := (x_1, \ldots, x_n)$ with $y(x_i) = y_i$. Then, for any $f$ defined on $\mathbf{x}$,

$$\|y - f\|_n^2 := \frac{1}{n} \sum_{i=1}^n |y_i - f(x_i)|^2 \tag{3.7}$$

is the empirical risk for $f$.

In order to estimate $f_\rho$ from the given data $\mathbf{z}$, we shall use the greedy algorithms OGA and RGA described in the previous section. We choose an arbitrary value of $a \geq 1$ and then fix it. We consider a dictionary $\mathcal{D}$ and truncations of this dictionary $\mathcal{D}_1, \mathcal{D}_2, \ldots$, as described in Section 2.3. We will use approximation from the span of the dictionary $\mathcal{D}_m$ in our algorithm, where we assume that the size $m$ is limited by

$$m \leq m(n) := \lfloor n^a \rfloor. \tag{3.8}$$

Our estimator is defined as follows:

(i) Given a data set **z** of size $n$, we apply the OGA, SPA or the RGA for the dictionary $\mathcal{D}_m$ to the function $y$ using the empirical inner product associated with the norm $\|\cdot\|_n$. In the case of the RGA, we use the second choice (2.6) for the parameter $\alpha_k$. This gives a sequence of functions $(\hat{f}_k)_{k=0}^\infty$ defined on **x**.

(ii) We define the estimator $\hat{f} := T\hat{f}_{k^*}$, where $Tu := T_B \min\{B, |u|\}\mathrm{sgn}(u)$ is the truncation operator at level $B$ and $k^*$ is chosen to minimize (over all $k > 0$) the penalized empirical risk

$$(3.9) \qquad \|y - T\hat{f}_k\|_n^2 + \kappa \frac{k \log n}{n},$$

with $\kappa > 0$ a constant to be fixed later.

We shall make some remarks about this algorithm. First, note that for $k = 0$, the penalized risk is bounded by $B^2$ since $\hat{f}_0 = 0$ and $|y| \leq B$. This means that we need not run the greedy algorithm for values of $k$ larger than $Bn/\kappa$. Second, our notation $\hat{f}$ suppresses the dependence of the estimator on **z**, which is again customary notation in statistics. The application of the $k$th step of the greedy algorithms requires the evaluation of $O(n^a)$ inner products. In the case of the OGA, we also need to compute the projection of $y$ onto a $k$-dimensional space. This could be done by doing Gram–Schmidt orthogonalization. Assuming that we had already computed an orthonormal system for step $k-1$, this would require an additional evaluation of $k-1$ inner products and then a normalization step. Finally, the truncation of the dictionary $\mathcal{D}$ is not strictly needed in some more specific cases such as neural networks (see Section 4).

In the following, we want to analyze the performance of our algorithm. For this analysis, we need to assume something about $f_\rho$. To impose conditions on $f_\rho$, we shall also view the elements of the dictionary normalized in the $L_2(\rho_X)$ norm $\|\cdot\|$. With this normalization, we denote by $\mathcal{L}_1$, $\mathcal{B}_p$, $\mathcal{L}_{1,r}$ and $\mathcal{B}_{p,r}$ the spaces of functions that have been previously introduced for a general Hilbert spaces $\mathcal{H}$. Here, we have $\mathcal{H} = L_2(\rho_X)$.

Finally, we denote by $\mathcal{L}_1^n$ the space of functions which admit an $\ell_1$ expansion in the dictionary when the elements are normalized in the empirical norm $\|\cdot\|_n$. This space is again equipped with a norm defined as the smallest $\ell_1$ norm among every admissible expansion. Similarly to $\|\cdot\|_n$, this norm depends on the realization of the design **x**.

3.2. *Error analysis.* In this section, we establish our main result, which will allow us, in the next section, to analyze the performance of the estimator under various smoothness assumptions on $f_\rho$.

THEOREM 3.1. *There exists $\kappa_0$ depending only on $B$ and $a$ such that if $\kappa \geq \kappa_0$, then for all $k > 0$ and for all functions $h$ in $\mathrm{Span}(\mathcal{D}_m)$, the estimator*

*satisfies*

$$(3.10) \qquad E(\|\hat{f} - f_\rho\|^2) \leq 8\frac{\|h\|^2_{\mathcal{L}_1(\mathcal{D}_m)}}{k} + 2\|f_\rho - h\|^2 + C\frac{k\log n}{n},$$

*where the constant $C$ only depends on $\kappa$, $B$ and $a$.*

The proof of Theorem 3.1 relies on a few preliminary results that we collect below. The first is a direct application of Theorem 3 from [18] or Theorem 11.4 from [11].

LEMMA 3.2. *Let $\mathcal{F}$ be a class of functions which are all bounded by $B$. For all $n$ and $\alpha, \beta > 0$, we have*

$$(3.11) \qquad \begin{aligned} &\Pr\{\exists f \in \mathcal{F} : \|f - f_\rho\|^2 \geq 2(\|y - f\|_n^2 - \|y - f_\rho\|_n^2) + \alpha + \beta\} \\ &\leq 14\sup_{\mathbf{x}} \mathcal{N}\left(\frac{\beta}{40B}, \mathcal{F}, L_1(\nu_\mathbf{x})\right) \exp\left(-\frac{\alpha n}{2568B^4}\right), \end{aligned}$$

*where $\mathbf{x} = (x_1, \ldots, x_n) \in X^n$ and $\mathcal{N}(t, \mathcal{F}, L_1(\nu_\mathbf{x}))$ is the covering number for the class $\mathcal{F}$ by balls of radius $t$ in $L_1(\nu_\mathbf{x})$, with $\nu_\mathbf{x} := \frac{1}{n}\sum_{i=1}^n \delta_{x_i}$ the empirical discrete measure.*

PROOF. This follows from Theorem 11.4 of [11] by taking $\epsilon = 1/2$ in that theorem. $\square$

We shall apply this result in the following setting. Given any set $\Lambda \subset \mathcal{D}$, we define $\mathcal{G}_\Lambda := \text{span}\{g : g \in \Lambda\}$ and denote by $T\mathcal{G}_\Lambda := \{Tf : f \in \mathcal{G}_\Lambda\}$ the set of all truncations of the elements of $\mathcal{G}_\Lambda$, where $T = T_B$ as before. We then define

$$(3.12) \qquad \mathcal{F}_k := \bigcup_{\Lambda \subset \mathcal{D}_m, \#(\Lambda) \leq k} T\mathcal{G}_\Lambda.$$

The following result gives an upper bound for the entropy numbers $\mathcal{N}(t, \mathcal{F}_k, L_1(\nu_\mathbf{x}))$.

LEMMA 3.3. *For any probability measure $\nu$, for any $t > 0$ and for any $\Lambda$ with cardinality $k$, we have the bound*

$$(3.13) \qquad \mathcal{N}(t, T\mathcal{G}_\Lambda, L_1(\nu)) \leq 3\left(\frac{2eB}{t}\log\frac{3eB}{t}\right)^{k+1}.$$

*Additionally,*

$$(3.14) \qquad \mathcal{N}(t, \mathcal{F}_k, L_1(\nu)) \leq 3n^{ak}\left(\frac{2eB}{t}\log\frac{3eB}{t}\right)^{k+1}.$$

PROOF. For each $\Lambda$ with cardinality $k$, Theorem 9.4 in [11] gives

$$\mathcal{N}(t, T\mathcal{G}_\Lambda, L_1(\nu)) \leq 3\left(\frac{2eB}{t}\log\frac{3eB}{t}\right)^{V_\Lambda}, \tag{3.15}$$

with $V_\Lambda$ the $VC$ dimension of the set of all subgraphs of $T\mathcal{G}_\Lambda$. It is easily seen that $V_\Lambda$ is smaller than the $VC$ dimension of the set of all subgraphs of $\mathcal{G}_\Lambda$, which, by Theorem 9.5 in [11], is less than $k+1$. This establishes (3.13). Since there are less than $n^{ak}$ possible sets $\Lambda$, the result (3.14) follows by taking the union of the coverings for all $T\mathcal{G}_\Lambda$ as a covering for $\mathcal{F}_k$. $\square$

Finally, we will need a result that relates the $\mathcal{L}_1$ and $\mathcal{L}_1^n$ norms.

LEMMA 3.4. *Given any dictionary $\mathcal{D}$, for all functions $h$ in $\mathrm{Span}(\mathcal{D})$, we have*

$$E(\|h\|_{\mathcal{L}_1^n}^2) \leq \|h\|_{\mathcal{L}_1}^2. \tag{3.16}$$

PROOF. We normalize the elements of the dictionary in $\|\cdot\| = \|\cdot\|_{\mathcal{H}}$. Given any $h = \sum_{g \in \mathcal{D}} c_g g$ and any $\mathbf{z}$ of length $n$, we write

$$h = \sum_{g \in \mathcal{D}} c_g g = \sum_{g \in \mathcal{D}} c_g^n \frac{g}{\|g\|_n}, \tag{3.17}$$

with $c_g^n := c_g \|g\|_n$. We then observe that

$$E\left(\left(\sum_{g \in \mathcal{D}} |c_g^n|\right)^2\right) = \sum_{(g,g') \in \mathcal{D} \times \mathcal{D}} |c_g c_{g'}| E(\|g\|_n \|g'\|_n)$$

$$\leq \sum_{(g,g') \in \mathcal{D} \times \mathcal{D}} |c_g c_{g'}| (E(\|g\|_n^2) E(\|g'\|_n^2))^{1/2}$$

$$= \sum_{(g,g') \in \mathcal{D} \times \mathcal{D}} |c_g c_{g'}| (\|g\|^2 \|g'\|^2)^{1/2}$$

$$= \sum_{(g,g') \in \mathcal{D} \times \mathcal{D}} |c_g c_{g'}|$$

$$= \left(\sum_{g \in \mathcal{D}} |c_g|\right)^2.$$

The result follows by taking the infimum over all possible admissible $(c_g)$ and using the fact that

$$E\left(\inf\left[\sum_{g \in \mathcal{D}} |c_g^n|\right]^2\right) \leq \inf E\left(\left[\sum_{g \in \mathcal{D}} |c_g^n|\right]^2\right). \tag{3.18}$$

$\square$

PROOF OF THEOREM 3.1. We write

$$\|\hat{f} - f_\rho\|^2 = T_1 + T_2, \tag{3.19}$$

where

$$T_1 := \|\hat{f} - f_\rho\|^2 - 2\left(\|y - \hat{f}\|_n^2 - \|y - f_\rho\|_n^2 + \kappa \frac{k^* \log n}{n}\right) \tag{3.20}$$

and

$$T_2 := 2\left(\|y - \hat{f}\|_n^2 - \|y - f_\rho\|_n^2 + \kappa \frac{k^* \log n}{n}\right). \tag{3.21}$$

From the definition of our estimator, we know that for all $k > 0$, we have

$$T_2 \leq 2\left(\|y - \hat{f}_k\|_n^2 - \|y - f_\rho\|_n^2 + \kappa \frac{k \log n}{n}\right). \tag{3.22}$$

Therefore, for all $k > 0$ and $h \in L_2(\rho_X)$, we have $T_2 \leq T_3 + T_4$ with

$$T_3 := 2(\|y - \hat{f}_k\|_n^2 - \|y - h\|_n^2) \tag{3.23}$$

and

$$T_4 := 2(\|y - h\|_n^2 - \|y - f_\rho\|_n^2) + 2\kappa \frac{k \log n}{n}. \tag{3.24}$$

We now bound the expectations of $T_1$, $T_3$ and $T_4$. For the last term, we have

$$\begin{aligned} E(T_4) &= 2E(|y - h(x)|^2 - |y - f_\rho(x)|^2) + 2\kappa \frac{k \log n}{n} \\ &= 2\|f_\rho - h\|^2 + 2\kappa \frac{k \log n}{n}. \end{aligned} \tag{3.25}$$

For $T_3$, we know from Theorems 2.3 and 2.4 that we have

$$\|y - \hat{f}_k\|_n^2 - \|y - h\|_n^2 \leq 4 \frac{\|h\|_{\mathcal{L}_1^n}^2}{k}. \tag{3.26}$$

Using, in addition, Lemma 3.4, we thus obtain

$$E(T_3) \leq 8 \frac{\|h\|_{\mathcal{L}_1}^2}{k}. \tag{3.27}$$

For $T_1$, we introduce $\Omega$, the set of $\mathbf{z} \in Z^n$ for which

$$\|\hat{f} - f_\rho\|^2 \geq 2\left(\|y - \hat{f}\|_n^2 - \|y - f_\rho\|_n^2 + \kappa \frac{k^* \log n}{n}\right). \tag{3.28}$$

Since $T_1 \leq \|\hat{f} - f_\rho\|^2 + 2\|y - f_\rho\|_n^2 \leq 6B^2$, we have

$$E(T_1) \leq 6B^2 \Pr(\Omega). \tag{3.29}$$

We thus obtain that for all $k > 0$ and for all $h \in L_2(\rho_X)$, we have

$$(3.30) \quad E(\|\hat{f} - f_\rho\|^2) \leq 8\frac{\|h\|^2_{\mathcal{L}_1}}{k} + 2\|f_\rho - h\|^2 + 2\kappa\frac{k \log n}{n} + 6B^2 \Pr(\Omega).$$

It remains to bound $\Pr(\Omega)$. Since $k^*$ can take an arbitrary value depending on the sample realization, we simply control this quantity by the union bound

$$(3.31) \quad \sum_{1 \leq k \leq Bn/\kappa} \Pr\Big\{\exists f \in \mathcal{F}_k :$$

$$\|f - f_\rho\|^2 \geq 2(\|y - f\|_n^2 - \|y - f_\rho\|_n^2) + 2\kappa\frac{k \log n}{n}\Big\}.$$

Denoting by $p_k$ each term of this sum, we obtain, by Lemma 3.2, that

$$(3.32) \quad p_k \leq 14 \sup_{\mathbf{x}} \mathcal{N}\left(\frac{\beta}{40B}, \mathcal{F}_k, L_1(\nu_\mathbf{x})\right) \exp\left(-\frac{\alpha n}{2568 B^4}\right),$$

provided $\alpha + \beta \leq 2\kappa\frac{k \log n}{n}$. Assuming without loss of generality that $\kappa > 1$, we can take $\alpha := \kappa\frac{k \log n}{n}$ and $\beta = 1/n$, from which it follows that

$$(3.33) \quad p_k \leq 14 \sup_{\mathbf{x}} \mathcal{N}\left(\frac{1}{40Bn}, \mathcal{F}_k, L_1(\nu_\mathbf{x})\right) n^{-\kappa k/2568 B^4}.$$

Using Lemma 3.3, we finally obtain

$$(3.34) \quad p_k \leq C n^{ak} n^{2(k+1)} n^{-\kappa k / 2568 B^4},$$

so that by choosing $\kappa \geq \kappa_0 = 2568 B^4(a+5)$, we always have $p_k \leq Cn^{-2}$. It follows that

$$(3.35) \quad \Pr(\Omega) \leq \sum_{k \leq Bn/\kappa} p_k \leq \frac{C}{n}.$$

This contribution is therefore absorbed into the term $2\kappa\frac{k \log n}{n}$ in the main bound and this completes our proof. □

REMARK 3.5. The value $\kappa_0 = 2568 B^4(a+5)$ is a pessimistic estimate due to the large numerical constant. In practice, this may result in selecting too small a value for $k^*$. An alternative approach to the complexity penalty for choosing $k^*$ is the so-called "hold out" method, which, in the present case, would consist in (i) splitting the sample set into two independent subsets $\{1, \ldots, \tilde{n}\}$ and $\{\tilde{n}+1, \ldots, n\}$, (ii) using the first subset to build the sequence $(\hat{f}_k)_{k \geq 0}$ and (iii) using the second subset to select a proper value of $k^*$.

3.3. *Rates of convergence and universal consistency.* In this section, we apply Theorem 3.1 in several situations which correspond to different prior assumptions on $f_\rho$. In order to control the approximation error resulting from using the truncated dictionary $\mathcal{D}_m$ in our algorithm, we take an arbitrary but fixed $a > 0$ and take the size $m$ exactly of the order of $n^a$,

$$m = m(n) := \lfloor n^a \rfloor. \tag{3.36}$$

We first consider the case where $f_\rho \in \mathcal{L}_{1,r}$. In that case, we know that for all $m$, there exists $h \in \mathrm{Span}(\mathcal{D}_m)$ such that $\|h\|_{\mathcal{L}_1(\mathcal{D}_m)} \leq M$ and $\|f_\rho - h\| \leq Mm^{-r}$, with $M := \|f_\rho\|_{\mathcal{L}_{1,r}}$. Therefore, Theorem 3.1 yields

$$E(\|\hat{f} - f_\rho\|^2) \leq C \min_{k>0} \left( \frac{M^2}{k} + M^2 n^{-2ar} + \frac{k \log n}{n} \right). \tag{3.37}$$

In order that the effect of truncating the dictionary does not affect the estimation bound, we make the assumption that $2ar \geq 1$. This allows us to delete the middle term in (3.37). Note that this is not a strong additional restriction over $f_\rho \in \mathcal{L}_1$ since $a$ can be fixed arbitrarily large.

COROLLARY 3.6. *If $f_\rho \in \mathcal{L}_{1,r}$ with $r > 1/2a$, then*

$$E(\|\hat{f} - f_\rho\|^2) \leq C(1 + \|f_\rho\|_{\mathcal{L}_{1,r}}) \left( \frac{n}{\log n} \right)^{-1/2}. \tag{3.38}$$

PROOF. We take $k := \lceil (M+1)^2 \frac{n}{\log n} \rceil^{1/2}$ in (3.37) and obtain the desired result. □

We next consider the case where $f_\rho \in \mathcal{B}_{p,r}$. In that case, we know that for all $m$ and for all $t > 0$, there exists $h \in \mathrm{Span}(\mathcal{D}_m)$ such that $\|h\|_{\mathcal{L}_1} \leq Mt^{\theta-1}$ [see (2.62)] and $\|f_\rho - h\| \leq M(t^\theta + t^{\theta-1} m^{-r})$ [see (2.64)], with $1/p = (1+\theta)/2$ and with $M = \|f\|_{B_{p,r}}$. Taking $t = k^{-1/2}$ and applying Theorem 3.1, we obtain

$$\begin{aligned} E(\|\hat{f} - f_\rho\|^2) \\ \leq C \min_{k>0} \left( M^2 k^{-2s} + M^2 (k^{-s} + k^{-s+1/2} n^{-ar})^2 + \frac{k \log n}{n} \right), \end{aligned} \tag{3.39}$$

with $s = 1/p - 1/2$. We now impose that $ar \geq 1/2$, which allows us to delete the term involving $n^{-ar}$.

COROLLARY 3.7. *If $f_\rho \in \mathcal{B}_{p,r}$ with $r \geq 1/(2a)$, then*

$$\begin{aligned} E(\|\hat{f} - f_\rho\|^2) \\ \leq C(1 + \|f_\rho\|_{\mathcal{B}_{p,r}})^{2/(2s+1)} \left( \frac{n}{\log n} \right)^{-2s/(1+2s)}, \end{aligned} \tag{3.40}$$

*with $C$ a constant depending only on $\kappa$, $B$ and $a$.*

PROOF. We take $k := \lceil (M+1)^2 \frac{n}{\log n} \rceil^{1/(1+2s)}$ in (3.39) and obtain the desired result. □

Let us finally prove that the estimator is universally consistent. For this, we need to assume that the dictionary $\mathcal{D}$ is complete in $L_2(\rho_X)$.

THEOREM 3.8. *For an arbitrary regression function, we have*

$$\lim_{n \to +\infty} E(\|\hat{f} - f_\rho\|^2) = 0. \tag{3.41}$$

PROOF. For any $\varepsilon > 0$ and $n$ sufficiently large, there exists $h \in \text{Span}(\mathcal{D}_m)$ which satisfies $\|f_\rho - h\| \leq \varepsilon$. According to Theorem 3.1, we thus have

$$E(\|\hat{f} - f_\rho\|^2) \leq C \min_{k>0} \left( \frac{\|h\|_{\mathcal{L}_1}^2}{k} + \varepsilon^2 + \frac{k \log n}{n} \right). \tag{3.42}$$

Taking $k = n^{1/2}$, this gives

$$E(\|\hat{f} - f_\rho\|^2) \leq C(\varepsilon^2 + n^{-1/2} \log n), \tag{3.43}$$

which is smaller than $2C\varepsilon^2$ for sufficiently large $n$. □

Finally, we remark that although our results for the learning problem are stated and proved for the OGA and RGA, they hold equally well when the SPA is employed.

**4. Neural networks.** Neural networks have been one of the main motivations for the use of greedy algorithms in statistics [2, 4, 14, 18]. They correspond to a particular type of dictionary. One begins with a univariate positive and increasing function $\sigma$ that satisfies $\sigma(-\infty) = 0$ and $\sigma(+\infty) = 1$ and defines the dictionary consisting of all functions

$$x \mapsto \sigma(\langle v, x \rangle + w) \tag{4.1}$$

for all vectors $v \in \mathbb{R}^D$ and scalars $w \in \mathbb{R}$, where $D$ is the dimension of the feature variable $x$. Typical choices for $\sigma$ are the Heaviside function $h = \chi_{x>0}$ or more general sigmoidal functions which are regularized versions of $h$.

In [18], the authors consider neural networks where $\sigma$ is the Heaviside function and the vectors $v$ are restricted to have at most $d$ nonzero coordinates ($d$-bounded fan-in) for some fixed $d \leq D$. We denote this dictionary by $\tilde{\mathcal{D}}$. With this choice of dictionary and using the standard relaxed greedy algorithm, they establish the convergence rate

$$E(\|\hat{f}_k - f_\rho\|^2 - \|f_\rho - f_a\|^2) \lesssim \frac{1}{k} + kd \frac{\log(Dn)}{n}, \tag{4.2}$$

where $f_a$ is the projection of $f_\rho$ onto the convex hull of the dictionary $\tilde{\mathcal{D}}$. This can also be expressed by

$$E(\|\hat{f}_k - f_a\|^2) \lesssim \frac{1}{k} + kd\frac{\log(Dn)}{n}, \tag{4.3}$$

which shows that with the choice $k = n^{1/2}$, the estimator converges to $f_a$ with rate $n^{-1/2}$ up to a logarithmic factor. In particular, the algorithm is not universally consistent since it only converges to the regression function when it belongs to this convex hull.

4.1. *Convergence results.* Let us apply our results to this particular setting. We first want to note that in this case, it is, from a theoretical point of view, not necessary to truncate the dictionary $\tilde{\mathcal{D}}$ into a finite dictionary in order to achieve our theoretical results. The truncation of dictionaries was used in the proof of Theorem 3.1 to bound the covering numbers of the sets $\mathcal{F}_k$ through the bound established in Lemma 3.3. In the specific case of $\tilde{\mathcal{D}}$, one can bound these covering numbers without truncation. Let us note that in this case,

$$\mathcal{F}_k := \bigcup_{\Lambda \subset \tilde{\mathcal{D}}, \#(\Lambda) \leq k} T\mathcal{G}_\Lambda, \tag{4.4}$$

where we no longer have the restriction that $\Lambda$ is in $\tilde{\mathcal{D}}_m$.

LEMMA 4.1. *For the dictionary $\tilde{\mathcal{D}}$, any probability measure $\nu$ of the type $\nu = \frac{1}{n}\sum_{i=1}^n \delta_{x_i}$ and any $k > 0$ and $t > 0$, we have the bound*

$$\mathcal{N}(t, \mathcal{F}_k, L_1(\nu)) \leq 3(n+1)^{k(D+1)}\left(\frac{2eB}{t}\log\frac{3eB}{t}\right)^{k+1}, \tag{4.5}$$

*where the sets $\mathcal{F}_k$ are defined as in* (4.4).

PROOF. As in the proof of Lemma 3.3, we first remark that

$$\mathcal{N}(t, T\mathcal{G}_\Lambda, L_1(\nu)) \leq 3\left(\frac{2eB}{t}\log\frac{3eB}{t}\right)^{k+1}. \tag{4.6}$$

We next use two facts from Vapnik–Chervonenkis theory (see, e.g., Chapter 9 in [11]): (i) if $\mathcal{A}$ is a collection of sets with VC dimension $\lambda$, then there are at most $(n+1)^\lambda$ sets of $\mathcal{A}$ separating the points $(x_1, \ldots, x_n)$ in different ways; (ii) the VC dimension of the collection of half-hyperplanes in $\mathbb{R}^D$ has VC dimension $D + 1$. It follows that there are at most $(n+1)^{D+1}$ hyperplanes separating the points $(x_1, \ldots, x_n)$ in different ways and therefore there are at most $(n+1)^{k(D+1)}$ ways of selecting $(g_1, \ldots, g_k)$ in $\mathcal{D}$ which will give different functions on the sample $(x_1, \ldots, x_n)$. The result follows by taking the union of the coverings on all possible $k$-dimensional subspaces. □

REMARK 4.2. Under the $d$-bounded fan-in assumption, the factor $(n+1)^{k(D+1)}$ can be reduced to $D^{kd}(\frac{en}{d+1})^{k(d+1)}$; see the proof of Lemma 3 in [18].

Based on this bound, a brief inspection of the proof of Theorem 3.1 shows that its conclusion still holds, now with $\kappa_0$ depending on $B$ and $D$. It follows that the rates of convergence in Corollaries 3.6 and 3.7 now hold under the sole assumptions that $f \in \mathcal{L}_1$ and $f \in \mathcal{B}_p$, respectively. On the other hand, the universal consistency result in Theorem 3.8 requires that the dictionary is complete in $L_2(\rho_X)$, which only holds when $d = D$, that is, when all direction vectors are considered.

Theorem 3.1 improves the bound (4.2) of [18] in two ways: on the one hand, $f_a$ is replaced by an arbitrary function $h$ which can be optimized and on the other hand, the value of $k$ can also be optimized.

REMARK 4.3. Note that truncating the dictionary is still necessary to obtain a practical scheme. Such a truncation can be achieved by restricting to a finite number $m$ of direction vectors $v$, which typically grows together with sample size $n$. In this case, we need to consider the spaces $\mathcal{L}_{1,r}$ and $\mathcal{B}_{p,r}$, which contain an additional smoothness assumption compared to $\mathcal{L}_1$ and $\mathcal{B}_p$. This additional amount of smoothness is meant to control the error resulting from the discretization of the direction vectors. We refer to [20] for general results on this problem.

4.2. *Smoothness conditions.* Finally, let us briefly discuss the meaning of the conditions $f \in \mathcal{L}_1$ and $f \in \mathcal{B}_p$ in the case of a dictionary $\mathcal{D}$ consisting of the functions (4.1) for a fixed $\sigma$ and for all $v \in \mathbb{R}^D$ and $w \in \mathbb{R}$. The following can be deduced from a classical result obtained in [4]: assuming that the marginal distribution $\rho_X$ is supported in a ball $B_r := \{|x| \leq r\}$, for any function $f$ having a converging Fourier representation

$$(4.7) \qquad f(x) = \int \mathcal{F}f(\omega) e^{i\langle \omega, x \rangle} \, d\omega,$$

the smoothness condition

$$(4.8) \qquad \int |\omega| |\mathcal{F}f(\omega)| \, d\omega < +\infty$$

ensures that

$$(4.9) \qquad \|f\|_{\mathcal{L}_1} \leq (2rC_f + |f(0)|) \leq 2rC_f + \|f\|_{L_\infty},$$

with $C_f := \int |\omega| |\mathcal{F}f(\omega)| \, d\omega$. Barron actually proves that condition (4.8) ensures that $f(x) - f(0)$ lies in the closure of the convex hull of $\mathcal{D}$ multiplied by $2rC_f$, the closure being taken in $L_2(\rho_X)$ and the elements of the dictionary being normalized in $L_\infty$. The bound (4.9) follows by remarking that the $L_\infty$ norm controls the $L_2(\rho_X)$ norm.

We can therefore obtain smoothness conditions which ensure that $f \in \mathcal{B}_p$ by interpolating between the condition $\omega \mathcal{F}f(\omega) \in L_1$ and $f \in L_2(\rho_X)$.

In the particular case where $\rho_X$ is continuous with respect to the Lebesgue measure, that is, $\rho_X(A) \leq c|A|$, we know that a sufficient condition to have $f \in L_2(\rho_X)$ is given by $\mathcal{F}f \in L_2$.

We can then rewrite the two conditions that we want to interpolate as $|\omega|^{-1}|\mathcal{F}f(\omega)| \in L_1(|\omega|^2\,d\omega)$ and $|\omega|^{-1}|\mathcal{F}f(\omega)| \in L_2(|\omega|^2\,d\omega)$. Therefore, by standard interpolation arguments, we obtain that a sufficient condition for a bounded function $f$ to be in $\mathcal{B}_p$ is given by

$$(4.10) \qquad |\omega|^{-1}|\mathcal{F}f(\omega)| \in wL_p(|\omega|^2\,d\omega)$$

or, in other words,

$$(4.11) \qquad \sup_{\eta>0} \eta^p \int \chi_{\{|\mathcal{F}f(\omega)| \geq \eta|\omega|\}} |\omega|^2\,d\omega < +\infty.$$

A slightly stronger, but simpler, condition is

$$(4.12) \qquad |\omega|^{-1}|\mathcal{F}f(\omega)| \in L_p(|\omega|^2\,d\omega),$$

which reads

$$(4.13) \qquad \int |\omega|^{2-p}|\mathcal{F}f(\omega)|^p\,d\omega < +\infty.$$

When $\rho_X$ is arbitrary, a sufficient condition for $f \in L_2(\rho_X)$ is $\mathcal{F}f \in L_1$, which actually also ensures that $f \in L_\infty$. In that case, we can again apply standard interpolation arguments and obtain that a sufficient condition for a bounded function $f$ to be in $\mathcal{B}_p$ is given by

$$(4.14) \qquad \sup_{A>0} A^{1-2/p} \int_{|\omega| \geq A} |\mathcal{F}f(\omega)|\,d\omega < +\infty.$$

A slightly stronger, but simpler, condition is

$$(4.15) \qquad \int |\omega|^{2/p-1}|\mathcal{F}f(\omega)|\,d\omega < +\infty.$$

**Acknowledgments.** The authors wish to thank Professor Rich Baraniuk and Rice University who hosted two of us (A. C. and R. D.) during the fall of 2005, when much of this research was completed.

A. BARRON  
DEPARTMENT OF STATISTICS  
YALE UNIVERSITY  
P.O. BOX 208290  
NEW HAVEN, CONNECTICUT 06520-8290  
USA  
E-MAIL: andrew.barron@yale.edu  

A. COHEN  
LABORATOIRE JACQUES-LOUIS LIONS  
UNIVERSITÉ PIERRE ET MARIE CURIE 175  
RUE DU CHEVALERET  
75013 PARIS  
FRANCE  
E-MAIL: cohen@ann.jussieu.fr  

W. DAHMEN  
INSTITUT FÜR GEOMETRIE  
UND PRAKTISCHE MATHEMATIK  
RWTH AACHEN  
TEMPLERGRABEN 55  
D-52056 AACHEN  
GERMANY  
E-MAIL: dahmen@igpm.rwth-aachen.de  

R. DEVORE  
INDUSTRIAL MATHEMATICS INSTITUTE  
UNIVERSITY OF SOUTH CAROLINA  
COLUMBIA, SOUTH CAROLINA 29208  
USA  
E-MAIL: devore@math.sc.edu